\documentclass[11pt]{article}
\usepackage{caption}
\usepackage{epsf,epsfig,amsfonts,amsgen,amsmath,amstext,amsbsy,amsopn,amsthm}
\usepackage{ebezier,eepic}
\usepackage{color}
\usepackage{multirow}
\usepackage{tikz}
\usepackage{enumerate}
\usepackage{enumitem}
\usepackage{url}
\usepackage{pgf,tikz,pgfplots}
\usepackage{mathrsfs}
\usepackage{amssymb}
\usepackage[symbol]{footmisc}
\usetikzlibrary{arrows}
\usepackage{hyperref}
\hypersetup{hidelinks, colorlinks=true,allcolors=black,pdfstartview=Fit,breaklinks=true}

\newtheorem{theorem}{Theorem} [section]
\newtheorem{lemma}[theorem]{Lemma}
\newtheorem{corollary}[theorem]{Corollary}
\newtheorem{definition}[theorem]{Definition}

\newtheorem{conjecture}[theorem]{Conjecture}

\newtheorem{remark}[theorem]{Remark}
\newtheorem{proposition} [theorem]{Proposition}
\newtheorem{question}[theorem]{Question}

\def\vv{\boldsymbol{\mathbf \mu}}
\def\0{\emptyset}

\def\ex{\mathrm{ex}}

\newcommand{\1}{{\uppercase\expandafter{\romannumeral1}}}
\newcommand{\2}{{\uppercase\expandafter{\romannumeral2}}}

\let\svthefootnote\thefootnote
\newcommand\blankfootnote[1]{%
	\let\thefootnote\relax\footnotetext{#1}%
	\let\thefootnote\svthefootnote%
}

\begin{document}
\title{On local Tur\'an density problems of hypergraphs}

\author{Chunqiu Fang\thanks{School of Computer Science and Technology, Dongguan University of Technology, Dongguan, Guangdong 523808, China. Email: chunqiu@ustc.edu.cn.}
~~~~~
 Guorong Gao\thanks{School of Mathematical Sciences, University of Science and Technology of China,
Hefei, Anhui 230026, China. Email: gaoguorong@ustc.edu.cn.}
~~~~~
Jie Ma\thanks{School of Mathematical Sciences, University of Science and Technology of China,
Hefei, Anhui 230026, China. Email: jiema@ustc.edu.cn.}
~~~~~
Ge Song\thanks{School of Management, University of Science and Technology of China,
Hefei, Anhui 230026, China. Email: gsong@ustc.edu.cn.}}
\date{}

\maketitle

\begin{abstract}
For integers $q\ge p\ge r\ge2$, we say that an $r$-uniform hypergraph $H$ has property $(q,p)$, if for any $q$-vertex subset $Q$ of $V(H)$, there exists a $p$-vertex subset $P$ of $Q$ spanning a clique in $H$. Let $T_{r}(n,q,p)=\min\{ e(H): H\subset \binom{[n]}{r}, H \text{~has property~} (q,p)\}$. The local Tur\'an density about property $(q,p)$ in $r$-uniform hypergraphs is defined as $t_{r}(q,p)=\lim_{n\to \infty}T_{r}(n,q,p)/\binom{n}{r}$. Frankl, Huang and R\"odl [J. Comb. Theory, Ser. A, 177 (2021)] showed that $\lim_{p\to\infty}t_{r}(ap+1,p+1)=\frac{1}{a^{r-1}}$ for positive integer $a$ and  $t_{3}(2p+1,p+1)=\frac{1}{4}$ for all $p\ge 3$ and asked the question that determining the value of $\lim_{p\to\infty}t_{r}(\gamma p+1,p+1)$, where $\gamma\ge 1$ is a real number. Based on the study of hypergraph Tur\'an densities, we determine some exact values of local Tur\'an densities and answer their question partially; in particular, our results imply that the equality in their question about exact values does not hold in general.

\end{abstract}


\section{Introduction}

A hypergraph $H$ is a pair $(V, E)$, where $V$ is a finite set of elements, called vertices, and  $E$ is a set of nonempty subsets of $V$, called hyperedges. The hypergraph $H$ is said to be $r$-uniform if every hyperedge of $H$ is of size $r$. The simple graphs are the $2$-uniform hypergraphs. We denote the numbers of the vertices and the hyperedges of $H$ by $\nu(H)$ and $e(H)$, respectively. Let $\binom{V(H)}{r}$ be the family of all $r$-element subsets of $V(H)$. For an $r$-uniform hypergraph $H$, we use $\overline{H}$ to denote the complement of $H$, that is $V(\overline{H})=V(H)$ and $E(\overline{H})=\binom{V(H)}{r}\setminus E(H)$.  For a vertex set $U\subset V(H)$, denote by $H[U]$ the subhypergraph induced by $U$, that is $H[U]=\{e\in E(H): e\subset U\}.$ A set $I$ of vertices in a hypergraph $H$ is independent if $H[I]$ does not contain a hyperedge. The independence number of $H$, denoted by $\alpha(H)$, is the maximum cardinality of an independent set in $H$.  A vertex set $X$ in a hypergraph $H$ is a clique if $H[X]$ is a complete $r$-uniform hypergraph. The clique number of $H$ is the maximum cardinality of a clique in $H$. For two vertices $x,y \in V(H)$, we denote $d_{H}(x)=|\{e\in E(H): e\ni x\}|$ and $d_{H}(xy)=|\{e\in E(H): e \supset\{x,y\}\}|$. Throughout this paper, we write $[n]=\{1, 2, \ldots, n\}.$

\subsection{Tur\'an problems}

The study of Tur\'{a}n number is one of the central topics in extremal graph (hypergraph) theory. Let $\mathcal{F}$ be a family of $r$-uniform hypergraphs, we say that an $r$-uniform hypergraph $H$ is $\mathcal{F}$-free if $H$ contains no member of  $\mathcal{F}$ as a subhypergraph. The Tur\'an number of $\mathcal{F}$, denoted $\ex(n, \mathcal{F})$,  is the maximum number of hyperedges in an $\mathcal{F}$-free $r$-uniform hypergraph on $n$ vertices. We let $\pi(\mathcal{F}):=\lim_{n\to\infty}\frac{\ex(n, \mathcal{F})}{\binom{n}{r}}$ and call  $\pi(\mathcal{F})$ the Tur\'an density of $\mathcal{F}$.
Katona, Nemetz and Simonovits \cite{Katona} used an averaging argument to show that the Tur\'an density of any family of hypergraphs exists. If $\mathcal{F}=\{F\}$, we denote the Tur\'an number $\ex(n,\{F\})$  by $\ex(n,F)$ and Tur\'an density $\pi(\{F\})$ by $\pi(F)$, respectively.

The classic Tur\'an Theorem \cite{Turan} says that for $p\ge 2$, the Tur\'an number $\ex(n,K_{p+1})$ is uniquely attained by the complete balanced $p$-partite graph on $n$ vertices. The following celebrated  Erd\H os-Stone-Simonovits Theorem \cite{erdossim, ErdosStone} gives a tight estimate for the Tur\'an number of any family of graphs.
\begin{theorem} [Erd\H os-Stone-Simonovits, \cite{erdossim, ErdosStone}]\label{ESS}
For any family $\mathcal{F}$ of non-empty graphs, we have
$$\ex(n, \mathcal{F})=\left(\frac{p-1}{2p}+o(1)\right)n^{2},$$
where $p=\Psi(\mathcal{F})=\min\{\chi(F): F\in\mathcal{F}\}-1$ denotes the subchromatic number of $\mathcal{F}$.
\end{theorem}

Note that Theorem \ref{ESS} implies that $\pi(F)=\frac{\chi(F)-2}{\chi(F)-1}$ for any non-empty graph $F$. Extending Tur\'an's Theorem to $r$-uniform hypergraphs is one of the most challenging problems in extremal graph theory. Erd\H os offered a money prize for determining $\pi(K_{k}^{r})$ for at least one pair $r, k$ with $k>r\ge 3$. However, no Tur\'an density $\pi(K_{k}^{r})$ is known for any $k>r\ge 3$ yet.

\begin{conjecture} [Tur\'an, \cite{Turan}] \label{Turanconjecture}
For every integer $k\ge 4$,
\begin{align} \nonumber
\pi(K_{k}^{3})=1-\left(\frac{2}{k-1}\right)^{2}.
\end{align}
\end{conjecture}

The most famous case $k=4$ in Conjecture \ref{Turanconjecture}, which declares that $\pi(K_{4}^{3})=\frac{5}{9}$ has attracted a lot of interest and activity through the years. It is known that many constructions (e.g. see \cite{Brown0, Fon, Kostochka}) achieve the conjectured value. To date, Razborov \cite{Razborov} proved the best known upper bound $\pi(K_{4}^{3})\le 0.561666$ by using flag algebra.

\subsection{Local Tur\'an density}

We consider local Tur\'an densities which are defined as follows.

\begin{definition}[Frankl-Huang-R\"odl, \cite{Frankl2}]
For integers $q\ge p\ge r\ge 2$, we say that an $r$-uniform hypergraph $H=(V,E)$ has property $(q,p)$, if for every vertex set $Q\subset \binom{V}{q}$, there exists a vertex set  $P\subset \binom{Q}{p}$ spanning a clique in $H$, that is, $\binom{Q}{r}\subset H$.
Let $T_{r}(n,q,p)=\min\{e(H): H\subset \binom{[n]}{r},\, H\, \text{has property}\, (q,p)\}$.
Let $t_{r}(n,q,p)=\frac{T_{r}(n,q,p)}{\binom{n}{r}}$ and $t_{r}(q,p):=\lim_{n\to \infty} t_{r}(n,q,p).$
We  call $t_{r}(q,p)$ the local Tur\'an density of property $(q,p)$ in $r$-uniform hypergraphs.
\end{definition}

For $q\ge p\ge r\ge 2$, let $\mathcal{G}^{r}_{q,p}$ be the family of $r$-uniform hypergraphs $G$ with $\nu(G)=q$ and $\alpha(G)\le p-1$.
Note that an $r$-uniform hypergraph $H$ has property $(q,p)$ if and only if $\overline{H}$ is $\mathcal{G}^{r}_{q,p}$-free. Then we have $T_{r}(n,q,p)+\ex(n, \mathcal{G}^{r}_{q,p})=\binom{n}{r}$. Thus,
 \begin{align}\label{relationship}
 t_{r}(q,p)=1-\pi(\mathcal{G}^{r}_{q,p}).
 \end{align}
This also shows that the local Tur\'an density exists. By the definition of local Tur\'an density, for any integers $q>p>r\ge 2$,
$t_{r}(q,p+1)\ge t_{r}(q+1,p+1)\ge t_{r}(q,p)\ge t_{r}(q+1,p).$

In \cite{ErdosSpencer}, it was shown that for $r=2$ (the graph case), $t_{2}(q,p)=1/\left\lfloor\frac{q-1}{p-1}\right\rfloor$. For general $r$, Frankl and Stechkin \cite{Frankl3} showed that if $q\le\frac{r}{r-1}(p-1)$, then $t_{r}(q,p)=1$.
Frankl \cite{Frankl1} proved that  $\lim_{p\to\infty} t_{3}(2p+1,p+1)=\frac{1}{4}$.
Frankl, Huang and R\"odl \cite{Frankl2} generalize it to the $r$-uniform hypergraphs. They proved the following Theorem.

\begin{theorem}[Frankl-Huang-R\"odl, \cite{Frankl2}]\label{integera}
For integers $r\ge 2$ and $a\ge 2$,
\begin{align}\nonumber
\lim_{p\to\infty} t_{r}(ap+1,p+1)=\frac{1}{a^{r-1}}.
\end{align}
\end{theorem}

In the same paper, they \cite{Frankl2} obtained the exact value of  $t_{3}(2p+1,p+1)$ for all $p\ge3$.
\begin{theorem} [Frankl-Huang-R\"odl, \cite{Frankl2}] \label{p=3,a=2}
For every integer $p\ge 3$,
\begin{align} \nonumber
t_{3}(2p+1,p+1)=\frac{1}{4}.
\end{align}
\end{theorem}

In light of these results, they ask the the following question.
\begin{question} [Frankl-Huang-R\"odl, \cite{Frankl2}]  \label{FHRquestion}
Is it possibly true that for every positive real number $\gamma>1$,
$$\lim_{p\to\infty}t_{r}(\gamma p+1,p+1)=1-\min_{F\in \mathcal{F}}\pi(F)=\frac{1}{\lfloor\gamma\rfloor^{r-1}},$$
where $\mathcal{F}$ is the family of all the $r$-uniform hypergraphs satifying $\nu(F)\ge \gamma\alpha(F)$?
\end{question}

\subsection{Our results}

In this section we introduce our main results. Our first result Theorem \ref{limgamma} gives a weaker version of the first equality of Question \ref{FHRquestion},
while the second result Theorem \ref{exactgamma} is about some exact values of local Tur\'{a}n densities, which implies that the second equality in Question \ref{FHRquestion} for limit values does not hold in general (on the other hand, it does hold for several sub-intervals).

\begin{theorem} \label{limgamma}
Let $r\ge 3$ be an integer and $\gamma >1$ be a real number. We have
$$\lim_{p\to\infty} t_{r}(\lfloor\gamma p\rfloor+1,p+1)=1-\pi(\mathcal{F}^{r}_{\gamma}),$$
where $\mathcal{F}^{r}_{\gamma}=\{r\text{-uniform graphs~} F: \nu(F)>\gamma\alpha(F)\}$. Furthermore, for any  integer $r\ge 3$, there exist a real number $\gamma$ and a positive integer $p_{0}=p_{0}(r, \gamma)$, such that for any integer $p\ge p_{0}$, we have  $t_{r}(\lfloor\gamma p\rfloor+1,p+1)=1-\pi(\mathcal{F}^{r}_{\gamma})$.
\end{theorem}

\noindent {\bf Remark 1.}
From Theorem \ref{limgamma}, after replacing the hypergraph family $\mathcal{F}$ in Question \ref{FHRquestion} by $\mathcal{F}^{r}_{\gamma}$, if we can show $\pi(\mathcal{F}^{r}_{\gamma})=\min_{F\in \mathcal{F}^{r}_{\gamma}}\pi(F)$,\footnote{For a general hypergraph family $\mathcal{G}$, it is known that $\pi(\mathcal{G})=\min_{G\in \mathcal{G}}\pi(G)$ does not hold.} then the first equality of Question \ref{FHRquestion} holds.

\begin{theorem}\label{exactgamma}
The following results about local Tur\'an density hold.

\noindent(1)~For a fixed positive real number $\gamma$, there exists a positive integer $p_{0}=p_{0}(\gamma)$, such that for all integers $p\ge p_{0}$, we have the following
\begin{align}  \nonumber
t_{3}(\lfloor\gamma p\rfloor+1,p+1)=
 \begin{cases}
\frac{7}{9}, &\text{~if~} \, \frac{3}{2}\le\gamma<\frac{5}{3},\\
\frac{19}{27}, &\text{~if~} \, \frac{5}{3}\le\gamma<\frac{7}{4},\\
\frac{1}{4}, &\text{~if~} \, 2\le\gamma<\frac{7}{3},\\
\frac{1}{9}, &\text{~if~} \, 3\le\gamma<\frac{22}{7}.
 \end{cases}
 \end{align}

\noindent(2)~Let $\gamma$ be a positive real number with $\frac{4}{3}\le \gamma <\frac{7}{5}$, there exists a positive integer $p_{0}=p_{0}(\gamma)$ such that for all integers $p\ge p_{0}$, we have $t_{4}(\lfloor\gamma p\rfloor+1,p+1)=1-\frac{4!}{4^{4}}$.

\noindent(3)~For integers $r\ge 3$ and $p\ge r^{2}-r-1$, $t_{r}(rp+1,(r-1)p+1)=1-\frac{r!}{r^{r}}$.
\end{theorem}

\noindent {\bf Remark 2.}
We point out that strictly speaking the first equality in Question~\ref{FHRquestion} may not be true. For example, let $r=3$ and $\gamma=\frac{3}{2}$. Then the hypergraph family $\mathcal{F}$ in Question~\ref{FHRquestion} contains the single edge as its member. Thus $1-\min_{F\in \mathcal{F}}\pi(F)=1$. But Theorem \ref{exactgamma} shows that
$t_{3}(\left\lfloor\frac{3}{2}p\right\rfloor+1,p+1)=\frac{7}{9}$, which implies that the first equality of Question~\ref{FHRquestion} fails in this case.
In Theorem \ref{limgamma}, we modify the corresponding definition to be $\nu(F)> \gamma\alpha(F)$.

\medskip

\noindent {\bf Remark 3.}
 It can be seen from the proofs that the conclusions of Theorems \ref{limgamma} and \ref{exactgamma} still hold when replacing $\lfloor\gamma p\rfloor$ with $\lceil\gamma p\rceil$.

\medskip

The rest of the paper is organized as follows. In Section \ref{localdensitylimitation}, we prove Theoem \ref{limgamma}. In Section \ref{localdensityexactvalue}, we give some extremal constructions for local Tur\'an densities and prove Theorem \ref{exactgamma}. In Section \ref{Concludingremarks}, we conclude with some remarks and nature problems on this topic.


\section{Proof of Theorem \ref{limgamma}}\label{localdensitylimitation}

In this section, we will prove Theorem \ref{limgamma}. The method comes from the proof of Theorem \ref{integera} in \cite{Frankl2}. For the integer pair $(q,p)$ with $q\le \gamma p$, we let $e(q,p)=\gamma p-q$. Note that since $q\ge p$, we always have $e(q,p)\le \gamma q-q=(\gamma-1)q$. We also need the following definition and lemmas.

\begin{definition} [Frankl-Huang-R\"odl, \cite{Frankl2}]
For $H\subset \binom{X}{r}$, the vertex set $Z\subset X$ is a $(w,v)$-hole if $|Z|=w>\gamma v$ and the clique number of $H[Z]$ is $v$.
\end{definition}

\begin{lemma}[Frankl-Huang-R\"odl, \cite{Frankl2}] \label{smallhole}
Suppose $H\subset \binom{X}{r}$ has property $(q,p)$ and $Z$ is a $(w,v)$-hole of $H$ with $w<q$, then $H[X\setminus Z]$ has property $(q-w,p-v)$.

\end{lemma}

\begin{proof}
Take an arbitrary set $U\in \binom{X\setminus Z}{q-w}$, then $U\cup Z \in \binom{X}{q}$. Since $H$ has property  $(q,p)$, $H[U\cup Z]$  contains a clique of size $p$. Hence $H[U]$ contains a clique of size $p-v$. Then $H[X\setminus Z]$ has property $(q-w,p-v)$.
\end{proof}

\begin{lemma}\label{holedensity}
For every $r\ge 3$, $\gamma >1$ and  $\epsilon>0$, there exists $\ell_{0}=\ell_{0}(r,\epsilon,\gamma)$ such that the following holds for all $\ell\ge \ell_{0}$. Suppose an $r$-graph $H$ on vertex set $X$ has property $(q,p)$ for all pairs $(q,p)$ with $q\le \ell$, $p=\lceil\frac{q}{\gamma}\rceil$ (In other words, $H$ does not have a $(w,v)$-hole with $\ell\ge w>\gamma v$). Then for all $Y\in \binom{X}{l}$,
\begin{align} \nonumber
e(H[Y])\ge (1-\epsilon)(1-\pi(\mathcal{F}^{r}_{\gamma}))\binom{\ell}{r}.
\end{align}
\end{lemma}

\begin{proof}
 Since
$\lim_{\ell\to\infty}\frac{\ex(\ell,\mathcal{F}^{r}_{\gamma})}{\binom{\ell}{r}}=\pi(\mathcal{F}^{r}_{\gamma}),$
 there exists $\ell_{0}=\ell_{0}(r,\epsilon,\gamma)$ such that for all $\ell\ge \ell_{0}$,
 $\ex(\ell,\mathcal{F}^{r}_{\gamma})\le (\pi(\mathcal{F}^{r}_{\gamma})+\epsilon(1-\pi(\mathcal{F}^{r}_{\gamma})))\binom{\ell}{r}$.
For $Y\subset \binom{X}{\ell}$, $H[Y]$ has property $(q,p)$ for all pairs $(q,p)$ with $q\le \ell$, $p=\lceil\frac{q}{\gamma}\rceil$, which means that $\overline{H[Y]}$ is $\mathcal{F}^{r}_{\gamma}$-free.  Thus, we have
$$e(H[Y])=\binom{\ell}{r}-e(\overline{H[Y]})\ge \binom{\ell}{r}-\ex(\ell,\mathcal{F}^{r}_{\gamma})\ge (1-\epsilon)(1-\pi(\mathcal{F}^{r}_{\gamma}))\binom{\ell}{r}.$$
\end{proof}

Now we are ready to prove Theorem \ref{limgamma}.

\begin{proof} [\bf Proof of Theorem \ref{limgamma}]
 Let $H$ be a maximum $\mathcal{F}^{r}_{\gamma}$-free $r$-graph on $n$ vertices, then $\overline{H}$ has property $(\lfloor\gamma p\rfloor+1,p+1)$ (also each of the properties $(\lceil\gamma p\rceil+1,p+1)$, $(\lfloor\gamma p\rfloor,p)$, and $(\lceil\gamma p\rceil,p)$). Thus, $t_{r}(\lfloor\gamma p\rfloor+1,p+1)\le1-\pi(\mathcal{F}^{r}_{\gamma}).$

Now we focus on the lower bound.  Given $\epsilon>0$, let us fix a large  integer $\ell\ge \ell_{0}= \ell_{0}(r,\frac{\epsilon}{2},\gamma)$, where $ \ell_{0}(r,\frac{\epsilon}{2},\gamma)$ is obtained from Lemma \ref{holedensity}. Let $$\theta_{0}=\min\{w-\gamma v: w>\gamma v, w\le \ell, v\ge r-1\}.$$
And then fix a much larger integer $L\ge \frac{4\gamma\ell^{2}}{\theta_{0}}$. Consider a sufficiently large $r$-uniform hypergraph $H\subset\binom{[n]}{r}$ having property $(q,p)$, where $q=\lfloor\gamma L\rfloor$ and $p=L$.
Our aim is to find a subset $X\subset [n]$ with $\binom{|X|}{r}>(1-\frac{\epsilon}{2})\binom{n}{r}$ such that $H[X]$ has no $(w,v)$-hole with $w\le \ell$ and $v\ge r-1$.

We start with $H_{0}=H$ and define $H_{i}$ inductively. Let $q_{0}=q$, $p_{0}=p$ and $X_{0}=[n]$. Suppose that $H_{i}\subset \binom{X_{i}}{r}$ has property $(q_{i},p_{i})$ and it still has a $(w_{i}, v_{i})$-hole. Then we take such a $(w_{i}, v_{i})$-hole $Z_{i}\subset X_{i}$ and set
\begin{align} \nonumber
X_{i+1}=X_{i}\setminus Z_{i},\qquad H_{i+1}=H_{i}[X_{i+1}].
\end{align}
By Lemma \ref{smallhole}, $H_{i+1}$ has property $(q_{i}-w_{i},p_{i}-v_{i})$. Moreover,
$$e(q_{i}-w_{i},p_{i}-v_{i})=\gamma(p_{i}-v_{i})-(q_{i}-w_{i})=\left(\gamma p_{i}-q_{i}\right)-\left(\gamma v_{i}-w_{i}\right)\ge e(q_{i},p_{i})+\theta_{0}. $$
Set $q_{i+1}=q_{i}-w_{i}$, $p_{i+1}=p_{i}-v_{i}$ and repeat. At every step, we have
\begin{align} \nonumber
\gamma(r-1)\le \gamma v_{i}<|X_{i}|-|X_{i-1}|=w_{i}\le \ell.
\end{align}

Since $v_{i}\ge r-1$ for all $i$ and $p_0=p$, we have $i\le \frac{p}{r-1}$.  Suppose at step $i$, the hypergraph $H_{i}$ no longer contains a $(w,v)$-hole with $w\le \ell$. Then we choose a subset $Q$ of size $\ell$ of $V(H_{i})$ uniformly at random. By Lemma \ref{holedensity}, we have
$$\frac{e(H_{i})}{\binom{|X_{i}|}{r}}=\frac{\mathbb{E}[e(H_{i}[Q])]}{\binom{\ell}{r}}\ge(1-\frac{\epsilon}{2})(1-\pi(\mathcal{F}^{r}_{\gamma})).$$
 On the other hand, $|X_{i}|\ge n-i\ell\ge n-\frac{p\ell}{r-1}$. Thus, for sufficiently large $n$, $\binom{|X_{i}|}{r}\ge (1-\frac{\epsilon}{2})\binom{n}{r}$. Therefore,
$$e(H)\ge e(H_{i})\ge \big(1-\frac{\epsilon}{2}\big)^{2}(1-\pi(\mathcal{F}^{r}_{\gamma}))\binom{n}{r}\ge(1-\epsilon)(1-\pi(\mathcal{F}^{r}_{\gamma}))\binom{n}{r}. $$

Otherwise suppose this process continues to produce $(w,v)$-holes. Let $m$ be the first index such that $q_{m}< 2l$. Since $e(q_{m},p_{m})\le (\gamma-1)q_{m}$ and $e(q_{i},p_{i})$ strictly increases at least $\theta_{0}$ after each step, it follows that $m\le \frac{(\gamma-1)q_{m}}{\theta_{0}}$. Thus,
$$\lfloor\gamma L\rfloor=q_{0}=q_{m}+\sum_{i=0}^{m-1}w_{i}\le 2\ell+m\ell<2\ell+\frac{2(\gamma-1)\ell^{2}}{\theta_{0}}< \frac{4\gamma \ell^{2}}{\theta_{0}}, $$
contradicting that $L\ge \frac{4\gamma\ell^{2}}{\theta_{0}}$. We finish the proof of the first part of Theorem \ref{limgamma}.
The proof of the other part is given by Theorem \ref{exactgamma}.
\end{proof}


\section{Exact values of some local Tur\'an densities}\label{localdensityexactvalue}

In this section, we determine some exact values of local Tur\'an densities based on the known hypergraph Tur\'an densities.

\subsection{Extremal constructions}\label{localdensityconstruction}

In this section, we construct two hypergraph families, which provide upper bounds for local Tur\'an densities.


\begin{definition}

For integers  $n, a\ge 1, r\ge 3$ and $k\ge 1$, let $\mathcal{K}^{r}_{n,a,k}$ be the family of $r$-uniform hypergraphs $H$ with $\nu(H)=n$ and $V(H)$ admits a vertex partition
$V(H)=V_{1}\cup V_{2}\cup \cdots \cup V_{a}$ and $V_{a}= U_{0}\cup U_{1}\cup \cdots \cup U_{k}$, such that
$E(H)=\left( \bigcup_{i=1}^{a-1}\binom{V_{i}}{r}\right)\bigcup\binom{V_{a}\setminus U_{0}}{r}\bigcup\left(\bigcup^{k}_{j=1}\binom{ U_{0}\cup U_{j}}{r}\right).$
Let $\rho_{r}(n,a,k)=\min\{e(H): H\in \mathcal{K}^{r}_{n,a,k}\}$ and $\rho_{r}(a,k)=\lim_{n\to\infty}\frac{\rho_{r}(n,a,k)}{\binom{n}{r}}$.
\end{definition}

When $k=1$, it is the extremal construction considered in \cite{Frankl2} and we have $\rho_{r}(a,1)=\frac{1}{a^{r-1}}.$
We compute its exact value for the case $r=3, a=1$ by the following proposition.

\begin{proposition} \label{rhovalue}
Let $k$ be an integer. We have $\rho_{3}(1,k)=\frac{5k+4}{9k}.$
\end{proposition}

\begin{proof}
Let $H\in \mathcal{K}^{3}_{n,1,k}$. Let $V(H)=U_{0}\cup U_{1}\cup \cdots \cup U_{k}$ and $E(H)=\binom{V(H)}{3}\setminus \bigcup_{1\le i_{1}<i_{2}\le k}U_{0}\times U_{i_{1}}\times U_{i_{2}}.$
We may assume that $|U_{1}|=|U_{2}|=\cdots=|U_{k}|=xn$. Then $|U_{0}|=(1-kx)n$ and $0\le x\le \frac{1}{k}$. Thus $e(H)=\binom{n}{3}-\binom{k}{2}(xn)^{2}(1-kx)n$. Denote $f(x)=x^{2}(1-kx)$ where $0\le x\le \frac{1}{k}$. By direct calculation, we have
$f_{\max}=f(\frac{2}{3k})=\frac{4}{27k^{2}}.$
Thus we have $\rho_{3}(1,k)=\lim_{n\to\infty}\frac{\binom{n}{3}-\binom{k}{2}f_{\max}n^{3}}{\binom{n}{3}}=\frac{5k+4}{9k}.$
\end{proof}




\begin{theorem} \label{upperbounds-rho}
Let $a, k\ge 1$, $r\ge3$ and $p\ge r-1$ be integers.  For large $n$ and any $r$-uniform  hypergraph $H\in \mathcal{K}^{r}_{n,a,k}$, $H$ has property $(\lfloor(a+1-\frac{1}{k})p\rfloor+1,p+1)$. Furthermore, we have $$t_{r}(\lfloor(a+1-\frac{1}{k})p\rfloor+1,p+1)\le \rho_{r}(a,k).$$
\end{theorem}

\begin{proof}
Let $V(H)=V_{1}\cup V_{2}\cup \cdots \cup V_{a} $ be the vertex set partition of $V(H)$, where $V_{a}=U_{0}\cup U_{1}\cup \cdots \cup U_{k}$. For any vertex set $X\in \binom{V(H)}{\lfloor(a+1-\frac{1}{k})p\rfloor+1}$, we will show that $H[X]$ contains a clique of size $p+1$.

We may assume that  $|X\cap V_{i}|\le p$ for any $1\le i\le a-1$, otherwise $H[X\cap V_{i}]$ contains a clique of size at least $p+1$. Thus, $|X\cap V_{a}|\ge 2p+1-\lfloor\frac{p}{k}\rfloor$.  If $p< 2k$, then $|X\cap V_{a}|\ge 2p$, one can easily find  a clique of size at least $p+1$ in $H[X\cap V_{a}]$. Now we assume that $p\ge2k$. Simiarly, we may assume that $p+1-\lfloor\frac{p}{k}\rfloor\le |X\cap U_{0}|\le p$, otherwise either $H[X\cap U_{0}]$ or $H[X\cap (U_{1}\cup \cdots \cup U_{k})]$ contains a clique of size at least $p+1$. Assume now that $|X\cap U_{0}|=p-\lfloor\frac{p}{k}\rfloor+t$, where $1\le t\le \lfloor\frac{p}{k}\rfloor$. Then  $|X\cap (U_{1}\cup \cdots \cup U_{k})|\ge p+1-t$. By the pigeonhole principle, there is some $j$ with  $1\le j\le k$ such that $|X\cap U_{j}|\ge \lceil\frac{p+1-t}{k}\rceil$. Thus, $|X\cap (U_{0}\cup U_{j})|\ge p-\lfloor\frac{p}{k}\rfloor+t+\lceil\frac{p+1-t}{k}\rceil\ge p+1$. Therefore $H[X\cap(U_{0}\cup U_{j})]$ contains a clique of size at least $p+1$.
\end{proof}



\begin{definition}

For integers  $n, a\ge 1, r\ge 3$ and $2\le k\le r-1$, let $\mathcal{L}^{r}_{n,a,k}$ be the family of $r$-uniform hypergraphs $H$ with $\nu(H)=n$ and $V(H)$ admits a vertex partition
$V(H)=V_{1}\cup V_{2}\cup \cdots \cup V_{a}$ and $V_{a}= U_{0}\cup U_{1}\cup \cdots \cup U_{k},$ such that
$E(H)=\left(\bigcup_{i=1}^{a-1}\binom{V_{i}}{r}\right)\bigcup\left(\bigcup^{k}_{j=0}\binom{V_{a}\setminus U_{j}}{r}\right).$
Let $\eta_{r}(n,a,k)=\min\{e(H): H\in \mathcal{L}^{r}_{n,a,k}\}$
and $\eta_{r}(a,k)=\lim_{n\to\infty}\frac{\eta_{r}(n,a,k)}{\binom{n}{r}}.$
\end{definition}

We compute its exact value for the case $a=1$ by the following proposition.

\begin{proposition} \label{etavalue}
Let $r\ge 3$ and $2\le k\le r-1$. We have $\eta_{r}(1,k)=\sum_{i=1}^{k}(-1)^{i+1}\binom{k+1}{k+1-i}(\frac{k+1-i}{k+1})^{r}.$
\end{proposition}

\begin{proof}
Let $H\in \mathcal{L}^{3}_{n,1,k}$ with the minmum number of hyperedges. Let $V(H)=U_{0}\cup U_{1}\cup \cdots \cup U_{k}$ and $E(H)=\bigcup^{k}_{j=0}\binom{V(H)\setminus U_{j}}{r}.$
We may assume that $|U_{0}|$=$|U_{1}|=|U_{2}|=\cdots=|U_{k}|=\frac{n}{k+1}$.
By the inclusion-exclusion principle, we have
\begin{align} \nonumber
e(H)&=|\bigcup^{k}_{j=0}E(H[V(H)\setminus U_{j}])|\\ \nonumber
&=\sum^{k}_{i=1}(-1)^{i+1}\sum_{0\le \ell_{1}<\ell_{2}<\cdots<\ell_{k+1-i}\le k}
e(H[U_{\ell_{1}}\cup U_{\ell_{2}}\cup\cdots\cup U_{\ell_{k+1-i}}])\\ \nonumber
&=\sum^{k}_{i=1}(-1)^{i+1}\binom{k+1}{k+1-i}\binom{\frac{k+1-i}{k+1}n}{r}.
\end{align}
Then the edge density of $H$ is $$\frac{\sum^{k}_{i=1}(-1)^{i+1}\binom{k+1}{k+1-i}\binom{\frac{k+1-i}{k+1}n}{r}}{\binom{n}{r}}=\sum_{i=1}^{k}(-1)^{i+1}\binom{k+1}{k+1-i}\left(\frac{k+1-i}{k+1}\right)^{r}+o(1), \text{~as~} n\to \infty.$$

\end{proof}

\begin{remark}\label{r/(r-1)}
Let $r\ge 3$. We have $\eta_{r}(1,r-1)=1-\frac{r!}{r^{r}}.$
\end{remark}




\begin{theorem} \label{upperbounds-eta}
Let $a\ge 1, r\ge 3$, $p\ge r-1$ and $2\le k\le r-1$ be integers. For large $n$ and any $r$-uniform  hypergraph $H\in \mathcal{L}^{r}_{n,a,k}$, $H$ has property $(\lfloor(a+\frac{1}{k})p\rfloor+1,p+1)$. Furthermore, we have $$t_{r}(\lfloor(a+\frac{1}{k})p\rfloor+1,p+1)\le \eta_{r}(a,k).$$
\end{theorem}

\begin{proof}
Let $V(H)=V_{1}\cup V_{2}\cup \cdots \cup V_{a} $ be the vertex set partition of $V(H)$, where $V_{a}=U_{0}\cup U_{1}\cup \cdots \cup U_{k}$. For any vertex set $X\in \binom{V(H)}{\lfloor(a+\frac{1}{k})p\rfloor+1}$, we will show that $H[X]$ contains a clique of size $p+1$.

We may assume that  $|X\cap V_{i}|\le p$ for any $1\le i\le a-1$, otherwise $X\cap V_{i}$ contains a clique of size at least $p+1$. Thus,  $|X\cap V_{a}|\ge \lfloor\frac{(k+1)p}{k}\rfloor+1\ge \frac{(k+1)p+1}{k}$.  By the pigeonhole principle, there is some $j$ with  $0\le j\le k$ such that $|X\cap (V_{a}\setminus U_{j})|\ge \lceil\frac{k}{k+1}|X\cap V_{a}|\rceil\ge p+1$. Thus  $H[X]$ contains a clique of size at least $p+1$.
\end{proof}


\subsection{Proof of Theorem \ref{exactgamma}}

In this subsection, we prove Theorem \ref{exactgamma}.
The key tool we use is the ``blow-up'' operation of the $r$-uniform hypergraphs. Let $H$ be an $r$-uniform hypergraph  on vertex set $\{v_{1}, v_{2}, \ldots, v_{l}\}$ and $\vv=(\mu_{1}, \mu_{2}, \ldots, \mu_{l})$ be an integer vector  with each $\mu_{i}\ge 1$.
Then the blow-up $H(\vv)$ is the $r$-uniform hypergraph formed by replacing the vertex $v_{i}$ of $H$ with a disjoint class of $\mu_{i}$ vertices for each $i$ and inserting a complete $r$-partite $r$-uniform hypergraph between any vertex classes corresponding to an edge in $H$. Given a family $\mathcal{H}=\{H_{1}, H_{2},\ldots, H_{s}\}$ of $r$-uniform hypergraphs and a family $\mathcal{T}=\{\vv^{1}, \vv^{2},\ldots, \vv^{s}\}$ of positive integer vectors with each $\vv^{i}$ of dimension $|V(H_{i})|$, we define the $\mathcal{T}$-blow-up of $\mathcal{H}$ to be $\mathcal{H}(\mathcal{T})=\{H_{i}(\vv^{i}): 1\le i\le s\}$. Brown and Simonovits \cite{Brown} proved the following extremely useful result.

\begin{theorem}[Brown-Simonovits, \cite{Brown}]\label{blowup}
If $\mathcal{H}=\{H_{1}, H_{2},\ldots, H_{s}\}$ is a family of $r$-uniform hypergraphs and $\mathcal{T}=\{\vv^{1}, \vv^{2},\ldots, \vv^{s}\}$ is a family of positive integer vectors with each $\vv^{i}$ of dimension $|V(H_{i})|$, then $\pi(\mathcal{H}(\mathcal{T}))=\pi(\mathcal{H})$.
\end{theorem}

We also need the following results on hypergraph Tur\'an densities.

\begin{theorem}[Frankl-F\"uredi, \cite{Frankl0}]\label{F5}
Let $F_{5}=\{123, 124, 345\}$, we have $\pi(F_{5})=\frac{2}{9}.$
\end{theorem}

\begin{theorem}[Baber-Talbot, \cite{Baber}]
Let $H_{1}=\{123, 124, 134, 234\}, H_{2}=\{123, 124, 125, 345, 346\},\\ H_{3}=\{123, 124, 345, 156, 256\}$, $H_{4}=\{123, 124, 125, 346, 356, 456\},$ then $\pi(\{H_{1}, H_{2}, H_{3}, H_{4}\})=\frac{8}{27}.$
\end{theorem}

Let $H^{-}_{2}=\{123, 124, 125, 345\}$.  Note that $H_{2}$ is a subgraph of some blow-up of $H^{-}_{2}$, thus we have the following corollary.

\begin{corollary}\label{H1234}
$\pi(\{H_{1}, H^{-}_{2}, H_{3}, H_{4}\})=\frac{8}{27}.$
\end{corollary}

\begin{theorem}[Mubayi-R\"odl, \cite{Mubayi}]\label{H7}
Let $H_{7}=\binom{[4]}{3}\cup\{(a,x,y) : a\in [4], x,y \in\{5,6,7\}, x\ne y\}\setminus\{(1,5,6)\}$, then $\pi(H_{7})=\frac{3}{4}.$
\end{theorem}

For $r\ge3$, let the generalized triangle $T^{r}$ be the $r$-uniform hypergraph on $(2r-1)$ vertices with three edges $\{1,\ldots, r\},\{1,\ldots, r-1, r+1\} \, \text{and} \, \{r, r+1,\ldots, 2r-1\}$. Pikhurko \cite{Pikhurko} obtained the Tur\'an density of $T^{4}$.

\begin{theorem}[Pikhurko, \cite{Pikhurko}]\label{T4}
$ \pi(T^{4})=\frac{4!}{4^{4}}.$
\end{theorem}

For $r\ge 3$, the generalized fan denoted by $F^{r}$ is the $r$-uniform hypergraph comprising $r+1$ edges $e_{1},\ldots,e_{r},e$ such that $e_{i}\cap e_{j}=\{x\}$ for all $i\ne j$, where $x\notin e$ and $|e_{i}\cap e|=1$ for all $i$. Mubayi and Pikhurko \cite{Mubayi0}  obtained the Tur\'an density of $F^{r}$.

\begin{theorem}[Mubayi-Pikhurko, \cite{Mubayi0}]\label{Fr}
$\pi(F^{r})=\frac{r!}{r^{r}}.$
\end{theorem}

Now we are ready to prove Theorem \ref{exactgamma}.

\begin{proof} [Proof of Theorem \ref{exactgamma}]
(1)~If $\frac{3}{2}\le\gamma<\frac{5}{3}$, by Theorem \ref{upperbounds-rho} and Proposition \ref{rhovalue},   $t_{3}(\lfloor\gamma p\rfloor+1,p+1)\le t_{3}(\lfloor\frac{3}{2}p\rfloor+1,p+1)\le\rho_{3}(1,2)=\frac{7}{9}.$
We just need to prove the lower bound. Let $\vv=(\frac{\lfloor\gamma p\rfloor+1}{5}, \frac{\lfloor\gamma p\rfloor+1}{5}, \frac{\lfloor\gamma p\rfloor+1}{5}, \frac{\lfloor\gamma p\rfloor+1}{5}, \frac{\lfloor\gamma p\rfloor+1}{5})$\footnote{When $\frac{\lfloor\gamma p\rfloor+1}{5}$ is not an integer, we always let it be its floor or ceil such that the sum of all the elements in $\vv$ equals $\lfloor\gamma p\rfloor+1$. We do the same if similar situations occur in the rest of the paper.}.
By Theorem \ref{blowup} and Theorem \ref{F5}, we have $\pi(F_{5}(\vv))=\pi(F_{5})=\frac{2}{9}$. Also, $\alpha(F_{5}(\vv))\le \frac{3(\lfloor\gamma p\rfloor+1)}{5}+O(1)<p+1$ for large $p$. By (\ref{relationship}), for large $p$, we have  $t_{3}(\lfloor\gamma p\rfloor+1,p+1)\ge 1-\pi(F_{5}(\vv))=\frac{7}{9}$.

If $\frac{5}{3}\le\gamma<\frac{7}{4}$, by Theorem \ref{upperbounds-rho} and Proposition \ref{rhovalue},  $t_{3}(\lfloor\gamma p\rfloor+1,p+1)\le t_{3}(\lfloor\frac{5}{3}p\rfloor+1,p+1)\le\rho_{3}(1,3)=\frac{19}{27}$. We just need to prove the lower bound. Let
\begin{align} \nonumber
\vv^{1}&=(\frac{\lfloor\gamma p\rfloor+1}{4}, \frac{\lfloor\gamma p\rfloor+1}{4}, \frac{\lfloor\gamma p\rfloor+1}{4}, \frac{\lfloor\gamma p\rfloor+1}{4}),\\ \nonumber
\vv^{2}&=(\frac{2(\lfloor\gamma p\rfloor+1)}{7}, \frac{2(\lfloor\gamma p\rfloor+1)}{7}, \frac{\lfloor\gamma p\rfloor+1}{7}, \frac{\lfloor\gamma p\rfloor+1}{7}, \frac{\lfloor\gamma p\rfloor+1}{7}),\\ \nonumber
\vv^{3}&=(\frac{\lfloor\gamma p\rfloor+1}{7}, \frac{\lfloor\gamma p\rfloor+1}{7}, \frac{\lfloor\gamma p\rfloor+1}{7}, \frac{\lfloor\gamma p\rfloor+1}{7}, \frac{2(\lfloor\gamma p\rfloor+1)}{7}, \frac{\lfloor\gamma p\rfloor+1}{7}),\\ \nonumber
\vv^{4}&=(\frac{\lfloor\gamma p\rfloor+1}{7}, \frac{\lfloor\gamma p\rfloor+1}{7}, \frac{\lfloor\gamma p\rfloor+1}{7}, \frac{\lfloor\gamma p\rfloor+1}{7}, \frac{\lfloor\gamma p\rfloor+1}{7}, \frac{2(\lfloor\gamma p\rfloor+1)}{7}). \nonumber
\end{align}
 Then for large $p$, $\alpha(H_{1}(\vv^{1}))\le\frac{\lfloor\gamma p\rfloor+1}{2}+O(1)<p+1$, $\alpha(H^{-}_{2}(\vv^{2}))\le\frac{4(\lfloor\gamma p\rfloor+1)}{7}+O(1)<p+1$, $\alpha(H_{3}(\vv^{3}))\le\frac{4(\lfloor\gamma p\rfloor+1)}{7}+O(1)<p+1$, $\alpha(H_{4}(\vv^{4}))\le\frac{4(\lfloor\gamma p\rfloor+1)}{7}+O(1)<p+1$.  By (\ref{relationship}), Theorem \ref{blowup} and Corollary \ref{H1234},  for large $p$, we have $t_{3}(\lfloor\gamma p\rfloor+1,p+1)\ge1-\pi(\{H_{1}(\vv^{1}), H^{-}_{2}(\vv^{2}), H_{3}(\vv^{3}), H_{4}(\vv^{4})\})=1-\frac{8}{27}=\frac{19}{27}$.

If $2\le\gamma<\frac{7}{3}$, by Theorem \ref{upperbounds-rho} and Proposition \ref{rhovalue},  $t_{3}(\lfloor\gamma p\rfloor+1,p+1)\le t_{3}(2p+1,p+1)\le \rho_{3}(2,1)=\frac{1}{4}$.
We just need to prove the lower bound. Let  $\vv=(\frac{\lfloor\gamma p\rfloor+1}{7}, \frac{\lfloor\gamma p\rfloor+1}{7}, \frac{\lfloor\gamma p\rfloor+1}{7}, \frac{\lfloor\gamma p\rfloor+1}{7}, \frac{\lfloor\gamma p\rfloor+1}{7}, \frac{\lfloor\gamma p\rfloor+1}{7}, \frac{\lfloor\gamma p\rfloor+1}{7})$. By Theorem \ref{blowup} and Theorem \ref{H7}, we have $\pi(H_{7}(\vv))=\pi(H_{7})=\frac{3}{4}.$ Also, $\alpha(H_{7}(\vv))\le\frac{3(\lfloor\gamma p\rfloor+1)}{7}+O(1)<p+1$ for large $p$.  By (\ref{relationship}), for large $p$,  we have $t_{3}(\lfloor\gamma p\rfloor+1,p+1)\ge1-\pi(H_{7}(\vv))=\frac{1}{4}$.

For $3\le \gamma<\frac{22}{7}$, the proof will be given in subsection \ref{specialTurandensity} (See Theorem \ref{gamma=3}).

(2)~For $\frac{4}{3}\le \gamma<\frac{7}{5}$, by Theorem \ref{upperbounds-eta} and Remark \ref{r/(r-1)}, $t_{4}(\lfloor\gamma p\rfloor+1,p+1)\le t_{4}(\lfloor\frac{4}{3}p\rfloor+1,p+1)\le \eta_{4}(1,3)=1-\frac{4!}{4^{4}}$. Now we focus on proving the lower bound.  Let $\vv=(\frac{\lfloor\gamma p\rfloor+1}{7}, \frac{\lfloor\gamma p\rfloor+1}{7}, \frac{\lfloor\gamma p\rfloor+1}{7},\frac{\lfloor\gamma p\rfloor+1}{7},\frac{\lfloor\gamma p\rfloor+1}{7},\frac{\lfloor\gamma p\rfloor+1}{7}, \frac{\lfloor\gamma p\rfloor+1}{7})$. By Theorem \ref{blowup} and Theorem \ref{T4}, we have $\pi(T^{4}(\vv))=\pi(T^{4})=\frac{4!}{4^{4}}$. Also, $\alpha(T^{4}(\vv))\le \frac{5(\lfloor\gamma p\rfloor+1)}{7}+O(1)<p+1$ for large $p$.  By (\ref{relationship}),  for large $p$,  we have $t_{4}(\lfloor\gamma p\rfloor+1,p+1)\ge 1-\pi(T^{4}(\vv))=1-\frac{4!}{4^{4}}$.

(3)~By Theorem \ref{upperbounds-eta} and Remark \ref{r/(r-1)}, $t_{r}(rp+1,(r-1)p+1)\le \eta_{r}(1,r-1)=1-\frac{r!}{r^{r}}$. Now we focus on proving the lower bound.
For $r\ge 3$ and $p\ge r^{2}-r-1$, let $F^{r}_p$ be the blow-up of $F^{r}$ by replacing each vertex of $e$ with $p-r^{2}+2r$ vertices, the common vertex $x$ with $(r-1)^{2}$ vertices and each vertex of $V(F^{r})\setminus(e\cup\{x\})$ with $r-1$ vertices.  By Theorem \ref{blowup} and Theorem \ref{Fr}, we have $\pi(F^{r}_p)=\pi(F^{r})=\frac{r!}{r^{r}}$. Note that $\nu(F^{r}_p)=rp+1$. Also, $\alpha(F^{r}_p)=(r-1)p$.  By (\ref{relationship}), we have  $t_{r}(rp+1,(r-1)p+1)\ge 1-\pi(F_{p}^{r})+O(1)=1-\frac{r!}{r^{r}}$.
\end{proof}

\subsection{Proof of Theorem \ref{exactgamma} when $3\le \gamma<\frac{22}{7}$} \label{specialTurandensity}

To obtain more exact values of local Tur\'an densities, we consider an extension of Theorem \ref{H7}.  We define the family $\mathcal{H}_{a}$ of $3$-uniform hypergraphs.

\begin{definition}
For an integer $a\ge 2$, let $P$ and $Q$ be two disjoint vertex sets with $|P|=2a$ and $Q=a+1$. Let $H_{a}$  be the $3$-uniform hypergraph on vertex set $P\cup Q$ with $E(H_{a})=\binom{P}{3}\cup\{xyz: x\in P, y,z\in Q\}$. Let
$$\mathcal{H}_{a}=\{H: \binom{P}{3} \subset E(H)\subset H_{a}, \text{~for all~} y, z\in Q, d_{H}(yz)\ge 2a-1\}.$$
\end{definition}

Using Tur\'an's Theorem \cite{Turan}, we have the following lemma about Tur\'an number in multigraphs.

\begin{lemma}\label{multigraph}
Let $n\ge a\ge 2$ be integers and $G$ be a multigraph on $n$ vertices. If any $a+1$ vertices of $G$ will induce at least one edge with multiplicities at least $2$, then $e(G)\ge\frac{2}{a}\binom{n}{2}-n.$
\end{lemma}

\begin{proof}
By Tur\'an's Theorem \cite{Turan}, $\ex(n,K_{a+1})\le(1-\frac{1}{a})\binom{n}{2}+\frac{n}{2}$. Thus $T_{2}(n,a+1,2)=\binom{n}{2}-\ex(n,K_{a+1})\ge\frac{1}{a}\binom{n}{2}-\frac{n}{2}$. Hence $e(G)\ge2\times T_{2}(n,a+1,2)\ge\frac{2}{a}\binom{n}{2}-n.$
\end{proof}

Now we will compute the Tur\'an density of family $\mathcal{H}_{a}$ under the condition that a warm version of Tur\'an's Conjecture (Conjecture \ref{Turanconjecture}) holds.

\begin{theorem}\label{Ha}
Let $a\ge2$ be an integer. If $\pi(K^{3}_{2a})\le 1-\frac{1}{a^{2}}$, then $\pi(\mathcal{H}_{a})=1-\frac{1}{a^{2}}$.
\end{theorem}

\begin{proof}
Let $G_{n,a}=\binom{[n]}{3}\setminus(\binom{V_{1}}{3}\cup\binom{V_{2}}{3}\cup\cdots\cup\binom{V_{a}}{3})$ with $V_{1}\sqcup V_{2}\sqcup\cdots\sqcup V_{a}=[n]$ and $V_{i}\in\{\lfloor\frac{n}{a}\rfloor, \lceil\frac{n}{a}\rceil\}.$ It is easy to check that $G_{n,a}$ is $\mathcal{H}_{a}$-free and $e(G_{n,a})=(1-\frac{1}{a^2}+o(1))\binom{n}{3}$ as $n\to\infty$. Thus $\pi(\mathcal{H}_{a})\ge1-\frac{1}{a^{2}}$.

For the upper bounds, we just need to show that for any real number $0<\varepsilon<\frac{1}{a^2}$, there exist an integer $N=N(\varepsilon)$ and a real number $C=C(N)$, such that for any integer $n\ge N$, $\ex(n,\mathcal{H}_{a})<(1-\frac{1}{a^{2}}+\varepsilon)\binom{n}{3}+an^{2}+C$. We prove it by induction on $n\ge N$. Since $\pi(K^{3}_{2a})\le 1-\frac{1}{a^{2}}$, we can choose a large integer $N_1=N_1(\varepsilon)$ such that $\ex(n,K^{3}_{2a})\le (1-\frac{1}{a^{2}}+\varepsilon)\binom{n}{3}$ holds for $n\ge N_1$. Furthermore, we can choose an integer $N=N(\varepsilon)\geq N_1$ and a real number $C$ such that
$(1-\frac{1}{a^{2}}+\varepsilon)\binom{N}{3}+aN^{2}+C=\binom{N}{3}$ and $(1-\frac{1}{a^{2}}+\varepsilon)\binom{n}{3}+an^{2}+C<\binom{n}{3}$ holds for $n>N$.

It is easy to see that our claim holds for $n=N$. Now let $n>N$ and assume that our claim holds for $n-1$. Let $G$ be a $3$-uniform hypergraph on vertex set $[n]$ with $e(G)\ge(1-\frac{1}{a^{2}}+\varepsilon)\binom{n}{3}+an^{2}+C$. Therefore, $G$ contains a clique with $2a$ vertices.  By symmetry, suppose $\binom{[2a]}{3}\subset G$. For $i\in[2a]$, we define the link graphs $G(i)=\left\{\{x,y\}\in\binom{[2a+1,n]}{2}: ixy\in G\right\}$. Let $M$ be the multigraph whose edge set is the union (with multiplicities) $G(1)\cup G(2)\cup\cdots\cup G(2a)$. Let $\overline{M}=\overline{G(1)}\cup \overline{G(2)}\cup\cdots\cup \overline{G(2a)}$. If there exists an $(a+1)$-subset $S\subset [2a+1,n]$, such that the multiplicity of each edge of  the induced subgraph of $\overline{M}$ on $S$ is no more than $1$, in other words, the multiplicity of each edge of  the induced subgraph of $M$ on $S$ is at least $2a-1$, then vertex set $[1,2a]\cup S$ induces  some copy of hypergraph in $\mathcal{H}_{a}$ in $G$.
Thus, we assume that the induced subgraph of $\overline{M}$ on any $(a+1)$-subset of $[2a+1,n]$ will contains at least one edge with multiplicities at least $2$. By Lemma \ref{multigraph}, $e(\overline{M})\ge\frac{2}{a}\binom{n-2a}{2}-(n-2a)$. Thus, $e(M)\le(2a-\frac{2}{a})\binom{n-2a}{2}+(n-2a)$. Therefore, there exists a vertex $i\in [2a]$, such that $e(G(i))\le\frac{1}{2a}e(M)\le(1-\frac{1}{a^2})\binom{n-2a}{2}+\frac{n}{2a}-1$.
It follows that
\begin{align}\nonumber
d_{G}(i)&\le\binom{2a-1}{2}+(2a-1)(n-2a)+(1-\frac{1}{a^2})\binom{n-2a}{2}+\frac{n}{2a}-1\\\nonumber
&\leq 2a^2+2an-4a^2+(1-\frac{1}{a^2})\binom{n-1}{2}\\\nonumber
&\le(1-\frac{1}{a^2}+\varepsilon)\binom{n-1}{2}+a(2n-1).\nonumber
\end{align}
By induction hypothesis, we have
\begin{align}\nonumber
e(G-i)=e(G)-d_{G}(i)\ge(1-\frac{1}{a^2}+\varepsilon)\binom{n-1}{3}+a(n-1)^{2}+C> \ex(n-1,\mathcal{H}_{a}).
\end{align}
So $G-i$ contains a copy of some hypergraph in $\mathcal{H}_{a}$. It follows that
$$\ex(n,\mathcal{H}_{a})<(1-\frac{1}{a^{2}}+\varepsilon)\binom{n}{3}+an^{2}+C.$$
\end{proof}

\begin{theorem}\label{a=3}
Let $a\ge2$ be an integer. If $\pi(\mathcal{H}_{a})=1-\frac{1}{a^{2}}$, then for any real number $\gamma$ with $a\le \gamma<a+\frac{1}{2a+1}$, there exist a postive integer $p_{0}=p(\gamma)$, such that for all integers $p\ge p_{0}$,
$$t_{3}(\lfloor\gamma p\rfloor+1,p+1)=\frac{1}{a^{2}}.$$
\end{theorem}

\begin{proof}
By Theorem \ref{upperbounds-rho}, $t_{3}(\lfloor\gamma p\rfloor+1,p+1)\le t_{3}(ap+1,p+1)\le\rho_{3}(a,1)=\frac{1}{a^{2}}$. We just need to prove the lower bound.  For each hypergraph $H\in \mathcal{H}_{a}$, let $H^{\lfloor\gamma p\rfloor+1}$ be the blow-up of $H$ by replacing each vertex of $P$ with $\frac{a(\lfloor\gamma p\rfloor+1)}{2a^{2}+a+1}$ vertices and replacing each vertex of $Q$ with $\frac{\lfloor\gamma p\rfloor+1}{2a^{2}+a+1}$ vertices. Let $\mathcal{H}^{\lfloor\gamma p\rfloor+1}_{a}=\{H^{\lfloor\gamma p\rfloor+1}: H\in \mathcal{H}_{a}\}$.  By Theorem \ref{blowup}, $\pi(\mathcal{H}^{\lfloor\gamma p\rfloor+1}_{a})=\pi(\mathcal{H}_{a})=1-\frac{1}{a^{2}}$.
For large $p$, $\alpha(H^{\lfloor\gamma p\rfloor+1})\le\frac{(2a+1)(\lfloor\gamma p\rfloor+1)}{2a^{2}+a+1}+O(1)<p+1$. By (\ref{relationship}), $t_{3}(\lfloor\gamma p\rfloor+1,p+1)\ge 1-\pi(\mathcal{H}^{\lfloor\gamma p\rfloor+1}_{a})=\frac{1}{a^{2}}.$
\end{proof}

Vaughan \cite{Vaughan} computed an upper bound of $\pi(K_{6}^{3})$ by  flag algebra.

\begin{theorem}[Vaughan, \cite{Vaughan}, \cite{Balogh}]\label{K6}
$\pi(K_{6}^{3})\le0.8583903.$
\end{theorem}

Combining Theorem \ref{Ha}, Theorem \ref{a=3} and Theorem \ref{K6}, we have the following conclusion.

\begin{theorem}\label{gamma=3}
Let $\gamma$ be a real number with $3\le \gamma <\frac{22}{7}$, there exist a postive integer $p_{0}=p(\gamma)$, such that for all integers $p\ge p_{0}$,
$$t_{3}(\lfloor\gamma p\rfloor+1,p+1)=\frac{1}{9}.$$
\end{theorem}

\begin{remark}
Combining Theorem \ref{Ha} and Theorem \ref{a=3}, if $\pi(K^{3}_{2a})\le 1-\frac{1}{a^{2}}$, then for any real number $\gamma$ with $a\le \gamma<a+\frac{1}{2a+1}$, there exist a postive integer $p_{0}=p(\gamma)$, such that for all integers $p\ge p_{0}$, $t_{3}(\lfloor\gamma p\rfloor+1,p+1)=\frac{1}{a^{2}}.$  If  Tur\'an's Conjecture (Conjecture \ref{Turanconjecture}) holds for odd number $k=2a+1$, that is $\pi(K^{3}_{2a+1})= 1-\frac{1}{a^{2}}$, then for any real number $\gamma$ with $a\le \gamma<a+\frac{1}{2}$, there exist a postive integer $p_{0}=p(\gamma)$, such that for all integers $p\ge p_{0}$, $t_{3}(\lfloor\gamma p\rfloor+1,p+1)=\frac{1}{a^{2}}.$

\end{remark}


\section{Concluding remarks}\label{Concludingremarks}

In this paper, we determine some exact values of local Tur\'an densities and answer Question \ref{FHRquestion} partially; in particular, our results imply that the equality in Question \ref{FHRquestion} about exact values of the limits does not hold in general.
The local Tur\'{a}n density problems are still widely open. Here we discuss some natural problems.

Frankl, Huang and R\"odl \cite{Frankl2} extended Theorem \ref{p=3,a=2} to $r$-uniform hypergraphs and made the following conjecture.

\begin{conjecture}[Frankl-Huang-R\"odl, \cite{Frankl2}]\label{FHRconj}
For integers $r\ge 2$ and $p$ sufficiently large,
\begin{align}\nonumber
t_{r}(2p+1,p+1)=\frac{1}{2^{r-1}}.
\end{align}

\end{conjecture}

Note that Question \ref{FHRquestion} is a generalization of Theorem \ref{integera} and closely related to Conjecture \ref{FHRconj}.
Based on Theorem \ref{limgamma}, we extend Question \ref{FHRquestion} to the following more general question. Recall that, for an integer $r\ge 3$ and a real number $\gamma >1$, we defined the family $\mathcal{F}^{r}_{\gamma}$ of $r$-uniform hypergraphs as
$\mathcal{F}^{r}_{\gamma}=\{r\text{-uniform graphs~} F: \nu(F)>\gamma\alpha(F)\}.$

\begin{question} \label{exactvaluegamma}
Let $r\ge 3$ be an integer and $\gamma >1$ be a real number.  Does there exist a positive integer $p_{0}=p_{0}(r,\gamma)$, such that for all integers $p\ge p_{0}$,
\begin{align}\nonumber
t_{r}(\gamma p+1,p+1)=1-\pi(\mathcal{F}^{r}_{\gamma})?
\end{align}
\end{question}

Using the blow-up theorem of Brown and Simonovits \cite{Brown}, one can show that Question \ref{exactvaluegamma} is equivalent to the following.

\begin{question}\label{compact}
Let $r\ge 3$ be an integer and $\gamma >1$ be a real number. Does there exist a finite subfamily $\mathcal{F}_{fin}\subset \mathcal{F}^{r}_{\gamma}$ such that
$$\pi(\mathcal{F}_{fin})=\pi(\mathcal{F}^{r}_{\gamma})?$$
\end{question}

We say a family $\mathcal{H}$ of $r$-uniform hypergraphs \emph{compact} if there is a finite subfamily $\mathcal{H}_{fin}\subset \mathcal{H}$ such that $\pi(\mathcal{H}_{fin})=\pi(\mathcal{H})$ (See Conjecture 1 in \cite{Brown}). In this notion, Question \ref{compact} asks whether the family $\mathcal{F}^{r}_{\gamma}$ is compact. The result of Frankl and Stechkin \cite{Frankl3} shows that the family $\mathcal{F}^{r}_{\gamma}$ is compact when $1\le\gamma<\frac{r}{r-1}$. Theorem \ref{p=3,a=2} shows that the family $\mathcal{F}^{3}_{2}$ is compact.
Moreover, the corresponding hypergraph family of each case considered in Theorem \ref{exactgamma} is also compact.

\bigskip

\noindent{\bf Acknowledgment.} Very recently we learn that Frankl and Nie \cite{FN} have also obtained results in this topic, which partially overlap with our Theorem 1.8 when $r=3$ and $2\leq \gamma<7/3$.

\bibliographystyle{unsrt}

\end{document}